\documentclass[11pt, reqno]{article}

\usepackage{amsmath}
\usepackage{latexsym}
\usepackage{amssymb}
\usepackage{amsfonts}
\usepackage{color}
\usepackage{arydshln}
\usepackage{graphics}
\usepackage{appendix}
\usepackage{setspace}

\newcommand{\proof}{\par\noindent{\it Proof.\ \ }}

\def\qed{\ifmmode\square\else\nolinebreak\hfill
$\square$\fi\par\vskip12pt}

\renewcommand{\proof}{\par\noindent{\it Proof.\ \ }}
\def\qed{\ifmmode\square\else\nolinebreak\hfill
$\square$\fi\par\vskip12pt}

\newtheorem{theorem}{Theorem}
\newtheorem{lemma}[theorem]{Lemma}%
\newtheorem{proposition}[theorem]{Proposition}%

\usepackage{authblk}

\title{\bf Construction of 2fi-optimal row-column designs}

\author{Yingnan Zhang}

\author{Jiangmin Pan}

\author{Lei Shi\thanks{Corresponding author. Email: lshi@ynu.edu.cn (Lei Shi)}}
\affil{
School of Statistics and Mathematics,
Yunnan University of Finance and Economics,
Kunming, China}

\date{}
\begin{document}
\maketitle

\noindent{\bf{Abstract:}}
Row-column factorial designs that provide unconfounded estimation of all main effects and the maximum number of  two-factor interactions (2fi's)
are called 2fi-optimal.  This issue has been paid great attention recently for its wide application in industrial or physical experiments. The constructions of
2fi-optimal two-level and three-level full factorial and fractional factorial row-column designs have been proposed.
However, the results for high prime level have not been achieved yet.
In this paper, we develop these constructions by giving a theoretical construction of $s^n$ full factorial
2fi-optimal row-column designs for any odd prime level $s$ and any parameter combination,
and theoretical constructions of $s^{n-1}$ fractional factorial
2fi-optimal row-column designs for any prime level $s$ and any parameter combination.

\section{Introduction}

Full and fractional factorial designs are widely used in industrial or physical experiments.
When the experimental units are not homogeneous, blocking is an effective method for reducing systematic variations
by arranging homogeneous experimental units into groups, and therefore improves  the precision of effect estimation.
For example, in industrial experiments, blocking factors can be site, batches of raw material and different days.
Factorial designs with only one blocking factor 
have been extensively investigated for many years.
The two blocking systems are referred to generally as row blocking and column blocking, and
the resulting designs are referred to as row-column designs.
Hinkelmann and Kempthorne (2005, Chapter 9) described  factorial row-column designs with batches of
raw material and machines as blocking factors. As an application,
Datta et al. (2017) described a veterinary trial
involving breeds and age groups of calves as blocking factors to study the effects of drugs.

The earliest investigation of row-column designs
for factorial experiments is by using Quasi-Latin squares introduced by Yates (1937).
For full and fractional factorial designs, a fundamental theoretical issue is how to choose good blocking
schemes and measure their ``goodness'' (Mukerjee and Wu, 1999).
The classic works in National Bureau of Standards (1957, 1959) contain many useful
blocking schemes for two-level and three-level full and fractional factorial designs.
Researchers have established  several  optimal criteria in different senses over the years.
For example, Cheng and Mukerjee (2003) concerned the construction of optimal row-column designs with
respect to maximum estimation capacity,
and Xu and Lau (2006) extended the criterion of minimum aberration (Fries and Hunter, 1980) to  blocked fractional factorial designs.

Row-column factorial designs which provide unconfounded estimation of all main effects and the maximum number of  two-factor interactions (2fi's)
are called 2fi-optimal.
Very recently, Zhou and Zhou (2023) 
found an upper bound of the number of unconfounded 2fi's,
and proposed a method for constructing 2fi-optimal row-column full $3^n$ factorial designs,
and remarked that the generalization of constructing 2fi-optimal design to higher prime level should be further studied.
Using this bound, some previous constructions for two-level full factorial row-column designs
in Choi and Gupta (2008), Wang (2017) and Godolphin (2019) can be proved to be 2fi-optimal.
For economic reasons, full factorial designs are rarely used in practice for more than seven factors
due to large number of runs.  Alternatively, fractional factorial designs
are the most widely used designs in experimentation investigations, refer to Suen et al. (1997).
For the fractional case,  Godolphin (2019) used trial and error to determine
which  design has good estimable properties.
She gave two-level and three-level fractional factorial designs with small runs,
based on the work of Chen et al. (1993),
and she indicated that ``the use of fractional factorials extends the level of complexity
substantially and only the surface of the problem is touched on here''.
Based on catalogue of Xu (2005, 2009), Zhou and Zhou (2023)  presented an computer algorithm
for constructing 2fi-optimal  row-column designs  for two-level and three-level fractional factorials.

The above contributions motivate a natural problem: construct 2fi-optimal row-column designs for high prime level $s$.
In this paper, we will present a theoretical construction of  $s^n$ full factorial
2fi-optimal row-column designs for any odd prime level $s$ and any parameter combination (in Section 3),
and theoretical constructions of $s^{n-1}$ fractional factorial
2fi-optimal row-column designs for any prime level $s$ and any parameter combination (in Section 4).
Some specific examples are given for illustrations in the Appendix.


\section{Preliminaries}

\subsection{Notation and definitions}\

We first introduce some notation  in finite projective geometry.
Let $s$ be a prime, and let $V_n$ be the $n$-dimensional vector space over $\rm{GF(s)}$,
the Galois field of order $s$.
Denote by $\rm{PG}(n-1,s)$ the set of all $1$-dimensional subspaces  of $V_{n}$,
which is called {\it projective geometry} of dimension $n-1$ over \rm{GF(s)}.
It is easily known that $\rm{PG}(n-1,s)$ contains exactly $(s^n-1)/(s-1)$ points (1-dimensional subspaces).
More information for finite projective geometry are referred to Street and Street (1987) or Raghavarao (1971).

An $s^n$ full factorial design has $n$ factors and $s$ levels with $s^n$ running comprising
all possible level combinations of the $n$ factors.
Factors are labelled by $F_1, F_2,\dots, F_n$, or alphabetically $A, B, \dots$ in examples.
An $s^{n-k}$ fractional factorial design is an $s^{-k}$ fraction of the $s^n$ design,
so it has $n$ factors but $s^{n-k}$ runs.
A word  (Suen et al. 1997) of an $s^{n-k}$ design is a $n$-dimensional vector with components in $\rm{GF(s)}$.
An $s^{n-k}$ fractional factorial design is determined by  $k$ linearly independent  words which generate
the group $\mathbb{Z}_s^k$, the elementary abelian group of order $s^k$.  This group is called  the treatment  defining contrast subgroup,
and accordingly, the $k$ words are called treatment defining words.
In the defining contrast subgroup, the word $(a_1,\dots,a_n)$ and all its multiples
$(\lambda a_1,\dots, \lambda a_n)$ with $\lambda\neq 0$ and $\lambda\in \rm{GF(s)}$
are considered to be the same.
That is, a word can be regarded as a point in $\rm{PG}(n-1,s)$. In this way,
the defining contrast  subgroup $\mathbb{Z}_s^k$ consists of $(s^k-1)/(s-1)$ nonzero distinct words.

The {\it resolution} (Box and Hunter, 1961a,b) of a fractional factorial design is the length of the shortest word in the
treatment defining contrast subgroup, excluding the identity $I$.
For example, the $2^{6-2}$ fractional factorial design with treatment defining words $ABCD$ and $CDEF$ has treatment defining
contrast subgroup $\{I, ABCD, CDEF, ABEF\}$, and the design has resolution IV.
For a design of resolution at least V, each main effect and two-factor interaction  are unconfounded with other main effects or two-factor interactions.
See Hinkelmann and Kempthorne (2005, Chapter 13)
for the differences  of resolution III, IV and V.

An $s^p\times s^q$ row-column design is a design in which level combinations of an $s^n$ factorial experiment
are arranged in an array with
$s^p$ rows and $s^q$ columns.
Each row or column of the array is regarded as a block.
The present paper focuses on $s^p\times s^q$ row-column designs without replication, i.e., $p+q\leq n$.
If $p+q=n$, the design is a full factorial row-column design, which will be discussed in Section 3.
If $p+q<n$, the design is an $s^{n-k}$ fractional factorial row-column design where $k=n-(p+q)$,
which will be investigated in Section 4.  Both cases have exactly one level combination
arranged in each cell (i.e., the intersection of a row and a column).

The interaction of $F_l$ and $F_m$ is denoted by $F_l\times F_m$, which can be partitioned into $s-1$ interaction
components $F_lF_m^v$, where $v=1,\dots, s-1$.
Let $Y_{ij}$ be the observation from the unit in the i-th row and j-th column. The analysis of variance model for factorial
row-column design is
\begin{equation}
Y_{ij}=\mu+\tau_i+\rho_j+\gamma_{d[i,j]}+\varepsilon_{ij},  \ \ \ \ \ i=1,\dots, s^p, \ \  j=1,\dots, s^q,
\end{equation}
where $\mu$ is the overall mean, $\tau_i$ is the effect of the i-th row, $\rho_j$ is the effect of the j-th column,
$d[i,j]$ is the level combination applied to unit $(i,j)$, and $\gamma_{d[i,j]}$ is the effect of $d[i,j]$.
The error terms $\varepsilon_{ij}$ are assumed to be uncorrelated, all with mean zero and variance $\sigma^2$.
Assume $d[i,j]=(x_1,\dots,x_n)$, and then the model for the treatment effect $\gamma_{d[i,j]}$ is

\begin{equation}
\gamma_{d[i,j]}=\gamma_{(x_1,\dots, x_n)}=\sum_{l=1}^n \alpha_{l,x_l}+\sum_{l=1}^n\sum_{m=2, m>l}^n \sum_{v=1}^{s-1} \beta_{lm, x_l+vx_m}^{(v)},
\end{equation}
where $\alpha_{l,u}$ is the effect of $F_l$ with level $u$,  $\beta_{lm,u}^{(v)}$ is the effect of the interaction component $F_lF_m^v$ with level $u$,
and $\sum_{u=0}^{s-1}\alpha_{l,u}=0$, $\sum_{u=0}^{s-1}\beta_{lm,u}^{(v)}=0$ for all $v=1,\dots, s-1$.
For example, if $s=3, n=4$ and $(x_1,x_2,x_3,x_4)=(0,1,1,2)$, then
$$
\begin{array}{rl}
\gamma_{(0,1,1,2)} =&\alpha_{1,0}+\alpha_{2,1}+\alpha_{3,1}+\alpha_{4,2}+\beta_{12,1}^{(1)}+\beta_{12,2}^{(2)}
+\beta_{13,1}^{(1)}+\beta_{13,2}^{(2)}+\beta_{14,2}^{(1)}+\beta_{14,1}^{(2)}+ \\
&\beta_{23,2}^{(1)}+\beta_{23,0}^{(2)}+ \beta_{24,0}^{(1)}+\beta_{24,2}^{(2)}+\beta_{34,0}^{(1)}+\beta_{34,2}^{(2)}.
\end{array}
$$

An effect is called {\it unconfounded} if its estimators are uncorrelated to any other estimators
in model (1). In particular, a main effect $F_l$ is called  {\it confounded} if the estimators of $\alpha_{l,u}$
with $u=0,\dots, s-1$
are correlated to some other estimators in model (1), and called unconfounded otherwise.
Similarly, an interaction component $F_l F_m^v$ is called unconfounded if the estimators of $\beta_{lm,u}^{(v)}$
with $u=0,\dots, s-1$ are uncorrelated to any other estimators in model (1).
An interaction $F_l\times F_m$ is unconfounded if all the $s-1$ components $F_lF_m^v$  are unconfounded.
This work concerns estimation of the main effects and two-factor interactions, and interactions of three or more factors are assumed negligible throughout this paper.

\subsection{Array generator matrix}\

An $s^p\times s^q$ row-column design with $n$ factors can be constructed from $\mathcal{G}$,
an array generator matrix of order  $(p+q)\times n$ with entries in $\rm{GF(s)}$,  by the following procedure.
The submatrices formed by the first $p$ rows and the last $q$ rows of $\mathcal{G}$ are denoted by $\mathcal{G}_c$ and  $\mathcal{G}_r$,
respectively.
Denote the rows of $\mathcal{G}_c$ by $\boldsymbol{u_i}=(u_{i,1},\dots, u_{i,n})$ with  $i=1,\dots, p$,
and the rows of $\mathcal{G}_r$ by $\boldsymbol{v_i}=(v_{i,1},\dots, v_{i,n})$ with $i=1,\dots, q$.
Assume $\mathcal{G}$ is full row rank.
Then $\boldsymbol{u_1},\dots,\boldsymbol{u_p}$ generate $\mathbb{Z}_s^p$, the $p$-dimensional vector space over $\rm{GF(s)}$,
and $\boldsymbol{v_1},\dots,\boldsymbol{v_q}$  generate $\mathbb{Z}_s^q$.
Set $\mathbb{Z}_s^p=\{\boldsymbol{x_i}|\boldsymbol{x_i}=(x_{i,1},\dots, x_{i,n}),  i=1,\dots, s^p\}$ with $\boldsymbol{x_1}=(0,\dots, 0)$,
and $\mathbb{Z}_s^q=\{\boldsymbol{y_j} |\boldsymbol{y_j}=(y_{j,1},\dots, y_{j,n}),  j=1,\dots, s^q\}$ with $\boldsymbol{y_1}=(0,\dots, 0)$.
Then $\mathcal{G}$ can be presented as

$$\mathcal{G}=\left(
\begin{array}{c}
\mathcal{G}_c\\  \hline
\mathcal{G}_r\\
\end{array}
\right)=\left(
\begin{array}{c}
\boldsymbol{u}_1\\
\vdots\\
\boldsymbol{u}_p\\ \hline
\boldsymbol{v}_1\\
\vdots\\
\boldsymbol{v}_q
\end{array}
\right)=
\left(
\begin{array}{cccc}
u_{1,1} & u_{1,2} &\dots &u_{1,n}\\
\vdots &&&\vdots\\
u_{p,1} & u_{p,2} &\dots &u_{p,n}\\ \hline
v_{1,1} & v_{1,2} &\dots &v_{1,n}\\
\vdots &&&\vdots\\
v_{q,1} & v_{q,2} &\dots &v_{q,n}
\end{array}
\right),$$
and the level combination in the $(i,j)$ unit of the row-column design is defined  to be $\boldsymbol{x_i}+\boldsymbol{y_j}$.
From an algebraic point of view,  each column block is a coset of $\mathbb{Z}_s^p$ in $\mathbb{Z}_s^{n-k}$,
and in particular the first column (column key block) of the design is $\mathbb{Z}_s^p$.
Similarly each row block is a coset of $\mathbb{Z}_s^q$ in $\mathbb{Z}_s^{n-k}$,
 with $\mathbb{Z}_s^q$ as the first row (row key block) of the design.
In fact, each level combination in the design is a linear combination of the rows of $\mathcal{G}$.
Therefore, $\mathcal{G}$ is required to be full row rank  to avoid replications.

\vskip0.1in
\textbf{Example 1.}
Let $s=3$, $p=3$, $q=2$, $n=7$,
and let $\mathcal{D}_1$  be the row-column design  constructed from  the following array generator matrix
$\mathcal{G}$:

$$\mathcal{G}=
\left(
\begin{array}{c}
\mathcal{G}_c\\  \hline
\mathcal{G}_r\\
\end{array}
\right)=\left(
\begin{array}{ccccccc}
1&0&0&2&2&1&0\\
1&1&2&1&2&0&0\\
2&2&2&2&0&0&1\\ \hline
1&1&1&0&1&0&1\\
0&1&2&1&0&1&1\\
\end{array}
\right).
$$
Then $\mathcal{D}_1$ is as in Table 1,
where rows and columns are randomized in  design $\mathcal{D}_1$.
The matrix $\mathcal{G}$ is full row rank, so there is no replication in $\mathcal{D}_1$.

\begin{table}
\caption{The $3^3\times 3^2$ row-column design $\mathcal{D}_1$ corresponding to $\mathcal{G}$ in Example 1.}
{\small
$$\begin{tabular}{c|cccccccc}
\hline
0000000 &1110101&0121011&1201112&1022120&2220202&0212022&2102221&2011210\\ \hline
1002210 &2112011&1120221& 2200022&2021000&0222112&1211202&0101101&0010120\\
1121200&2201001&1212211&2022012&2110020&0011102&1000222&0220121&0102110\\
2222001&0002102&2010012&0120110&0211121&1112200&2101020&1021222&1200211\\
2120110&0200211&2211121&0021222&0112200&1010012&2002102&1222001&1101020\\
0211010&1021111&0002021&1112122&1200100&2101212&0120002&2010201&2222220\\
0221211&1001012&0012222&1122020&1210001&2111110&0100200&2020102&2202121\\
2110212&0220010&2201220&0011021&0102002&1000111&2022201&1212100&1121122\\
0010210&1120011&0101221&1211022&1002000&2200112&0222202&2112101&2021120\\
2202202&0012000&2020210&0100011&0221022&1122101&2111221&1001120&1210112\\
2001120&0111221&2122101&0202202&0020210&1221022&2210112&1100011&1012000\\
2212100&0022201&2000111&0110212&0201220&1102002&2121122&1011021&1220010\\
1111002&2221100&1202010&2012111&2100122&0001201&1020021&0210220&0122212\\
1210220&2020021&1001201&2111002&2202010&0100122&1122212&0012111&0221100\\
0122020&1202121&0210001&1020102&1111110&2012222&0001012&2221211&2100200\\
0112122&1222220&0200100&1010201&1101212&2002021&0021111&2211010&2120002\\
1220121&2000222&1011102&2121200&2212211&0110020&1102110&0022012&0201001\\
0020120&1100221&0111101&1221202&1012210&2210022&0202112&2122011&2001000\\
1101101&2211202&1222112&2002210&2120221&0021000&1010120&0200022&0112011\\
1012111&2122212&1100122&2210220&2001201&0202010&1221100&0111002&0020021\\
2011021&0121122&2102002&0212100&0000111&1201220&2220010&1110212&1022201\\
2100011&0210112&2221022&0001120&0122101&1020210&2012000&1202202&1111221\\
0201112&1011210&0022120&1102221&1220202&2121011&0110101&2000000&2212022\\
2021222&0101020&2112200&0222001&0010012&1211121&2200211&1120110&1002102\\
1022012&2102110&1110020&2220121&2011102&0212211&1201001&0121200&0000222\\
1200022 &2010120&1021000&2101101&2222112&0120221&1112011&0002210&0211202\\
0102221&1212022&0220202&1000000&1121011&2022120&0011210&2201112&2110101\\

\hline
\end{tabular}$$
}
\end{table}

In order to find the resolution of design $\mathcal{D}_1$, we need to find two linearly independent treatment defining words.
In fact, for $s^{n-k}$ fractional factorial designs,
the level combinations that are used in designs are orthogonal to the treatment defining words.
Solving the following equations,
\begin{equation*}
\begin{cases}
x_1+2x_4+2x_5+x_6\equiv 0  \\
x_1+x_2+2x_3+x_4+2x_5\equiv 0  \\
2x_1+2x_2+2x_3+2x_4+x_7\equiv 0  \pmod 3\\
x_1+x_2+x_3+x_5+x_7 \equiv 0 \\
x_2+2x_3+x_4+x_6+x_7 \equiv 0
\end{cases}
\end{equation*}
we derive that two linearly
independent treatment defining words of $\mathcal{G}$ are $(0,1,1,1,1,2,0)$ and $(0,1,2,2,1,0,2)$,
which correspond to  interactions $BCDEF^2$ and $BC^2D^2EG^2$, respectively.
Therefore, all the  treatment defining words are $BCDEF^2$, $BC^2D^2EG^2$, $BEFG$ and $CDFG^2$,
which shows that the resolution of design $\mathcal{D}_1$ is IV.


\subsection{Identifying 2fi-optimal designs from array generator matrices}\

The problem for identifying confounded effects for the designs constructed by array generator matrices is
crucial for the construction of 2fi-optimal row-column designs.
The following lemma addressed this problem by using Cheng and Mukerjee (2003, Lemma 1),
where the $i$-th column of $\mathcal{G}$, $\mathcal{G}_c$ and $\mathcal{G}_r$
are denoted by $\mathcal{G}_i$, $\mathcal{G}_{c,i}$ and $\mathcal{G}_{r,i}$, respectively.

\begin{lemma}[Zhou and Zhou, 2023]\label{identify}
Let $\mathcal{G}$ be an array generator matrix. Let $\mathcal{D}$ be the row-column design constructed by $\mathcal{G}$.
Then the following statements hold.
\begin{itemize}
\item[(1)] A main effect or a two-factor interaction is confounded with another main effect or two-factor interaction if and
only if the corresponding columns of $\mathcal{G}$ of these factors are linearly dependent.
\item[(2)] If $\mathcal{G}_i\neq 0$, $F_i$ is confounded with column effects if and only if $\mathcal{G}_{c,i}=0$,
and it is confounded with row effects if and only if $\mathcal{G}_{r,i}=0$.
\item[(3)] If $\mathcal{G}_i$ and $\mathcal{G}_j$ are linearly independent, $F_i\times F_j$ is confounded with column effects if and only if
$\mathcal{G}_{c,i}$ and $\mathcal{G}_{c,j}$ are linearly dependent, and it is confounded with row effects if and only if
$\mathcal{G}_{r,i}$ and $\mathcal{G}_{r,j}$ are linearly dependent.
\end{itemize}
\end{lemma}

Lemma \ref{identify} provides an approach to identify confounded treatment effects.
Take  the matrix $\mathcal{G}$ given in Example 1 for illustration.
Since all treatment defining words in the example are $BCDEF^2, BC^2D^2EG^2, BEFG$ and  $CDFG^2$,
and only  $BEFG$ and $CDFG^2$ are of length $4$,
by Lemma \ref{identify}(1), the two-factor interactions which are unconfounded with other main effects or two-factor interactions are
$A\times B, A\times C, A\times D, A\times E$, $A\times F$, $A\times G$, $B\times C$, $B\times D$, $C\times E$ and $D\times E$.
Since the resolution of $\mathcal{D}_1$ is IV,  each  main  effect is not confound with other main effects or two-factor
interactions.
By Lemma \ref{identify} (2), it is straightforward to see that  no main effects
are confounded with row effects or column effects.
It then follows from Lemma \ref{identify}(3) that no two-factor interactions are confounded with column effects,
and $A\times E$, $B\times G$, $D\times F$ are confounded with row effects.
To sum up, for  design $\mathcal{D}_1$ in Example 1, all main effects are unconfounded,
and there are only $9$ unconfounded two-factor interactions,
which are $A\times B, A\times C, A\times D$, $A\times F$, $A\times G$, $B\times C$, $B\times D$, $C\times E$ and $D\times E$.

To evaluate the 2fi-optimality of an $s^p\times s^q$ row-column design with $s$ level and $n$ factors,
Zhou and Zhou (2023, Theorem 1) obtained an upper bound $\phi(s,p,q,n)$ of the number of unconfounded two-factor interactions:
$$\phi(s,p,q,n)=\tbinom{n}{ 2}-({s^{\min{(p,q)}}-1\over s-1}){\tbinom{\alpha}{2}}-\alpha\beta,$$
where $\alpha=\lfloor n/({s^{\min{(p,q)}}-1\over s-1})  \rfloor$ , $\beta=n-({s^{\min{(p,q)}}-1\over s-1})\alpha$,
and $\lfloor x\rfloor$ is the greatest integer not exceeding $x$.
For a design $\mathcal{D}$ with all main effects unconfounded,
they further defined  $t_{\mathcal{D}}/ \phi(s,p,q,n)$ as the 2fi-efficiency of the design  $\mathcal{D}$,
where $t_{\mathcal{D}}$ denotes the actual number of confounded two-factor interactions of $\mathcal{D}$.
The 2fi-efficiency may measure how far a  design is from 2fi-optimal.
Clearly, a design with 2fi-efficiency equal to $1$ is 2fi-optimal. But the converse is not true, that is,
a 2fi-optimal design does not guarantee that the 2fi-efficiency equals $1$, see
Zhou and Zhou (2023, Table 1) for some 2fi-optimal designs with 2fi-efficiency strictly less than $1$.
For illustration, we calculate the 2fi-efficiency in Example 1.
It is easy to see $\phi (3,3,2,7)= 18$, comparing to the total  $21$  two-factor interactions.
It has been observed above that design $\mathcal{D}_1$ has $9$ unconfounded two-factor interactions,
so design $\mathcal{D}_1$ has 2fi-efficiency $0.5$.

For designs with  resolution at least V,  the following proposition
gives  equivalent  conditions for identifying  which design has 2fi-efficiency $1$.

\begin{proposition}[Zhou and Zhou, 2023]\label{conditions}
For $p\geq 2$, let $\mathcal{G}$ be an array generator matrix such that any four columns are linearly independent if
$n\geq 4$, or $G$ has full column rank if $n\leq 3$. Matrix $\mathcal{G}$ yields a $2$fi-optimal row-column design whose
$2$fi-efficiency equals $1$ if and only if the following conditions are satisfied.
\vspace{0.05in}

{\it Condition $(1)$.}  
{ No column of $\mathcal{G}_c$ or $\mathcal{G}_r$ is a zero vector.}
\vspace{0.05in}

{\it Condition $(2)$.}  
{The number of columns of $\mathcal{G}_c$ belonging to each point of $\rm{PG(p-1,s)}$ is $\alpha$ or $\alpha+1$,
where $\alpha=\lfloor (s-1)n/(s^p-1)\rfloor$.}
\vspace{0.05in}

{\it Condition $(3)$.}  
{For $1\leq i<j\leq n$, $\mathcal{G}_{r,i}$ and $\mathcal{G}_{r,j}$ belong to different points of $\rm{PG(q-1,s)}$
if $\mathcal{G}_{c,i}$ and $\mathcal{G}_{c,j}$ belong to different points of $\rm{PG(p-1,s)}$.}
\end{proposition}

We next give a design with 2fi-efficiency $1$, which has
the same parameters as Example 1 except interchanging $p$ and $q$.

\vskip0.1in
\textbf{Example 2.}
Let $s=3$, $p=2$, $q=3$, $n=7$,
and let $\mathcal{D}_2$ be the row-column design  constructed from the array generator matrix $\mathcal{G}$:
$$\mathcal{G}=\left(
\begin{array}{ccccccc}
1&1&1&0&1&0&1\\
0&1&2&1&0&1&1\\ \hline
1&0&0&2&2&1&1\\
1&1&2&1&2&0&1\\
2&2&2&2&0&0&1\\
\end{array}
\right).
$$

The design $\mathcal{D}_2$ has resolution $V$, with the four defining words $BCDEF^2$,$AB^2D^2EG^2$,
$ACE^2F^2G^2$, and $ABC^2DFG^2$, i.e., any four columns of $\mathcal{G}$ are linearly independent,
so it satisfies the assumption of Proposition \ref{conditions}.
Apparently, no column of $\mathcal{G}_c$ or $\mathcal{G}_r$ is a zero vector.
There are four points of $\rm{PG(1,3)}$, denoted by $(1,0), (0,1), (1,1), (1,2)$,
which are the four $1$-dimensional subspaces of the $2$-dimensional
vector space over $\rm{GF(3)}$.
For Condition (2) of Proposition \ref{conditions}, $\alpha=\lfloor7/4\rfloor=1$,
and the numbers of columns that belong to each of the 4 points are $2,2,2,1$, respectively,
so Condition (2) is satisfied.
Since the columns of $\mathcal{G}_r$ represent different points of $\rm{PG(2,3)}$,
Condition (3) also holds. Therefore,  design $\mathcal{D}_2$ is 2fi-optimal with
2fi-efficiency $1$.

Proposition \ref{conditions}  is to identify  if a design achieves 2fi-efficiency  $1$ under the assumption $p\geq 2$.
For $p=1$, one easily derives a similar proposition from Lemma \ref{identify}.

\begin{proposition}\label{p=1}
For $p=1$, let $\mathcal{G}$ be an array generator matrix such that any three columns
are linearly independent if $n\geq 3$,
or $G$ has full column rank if $n=2$.
Matrix $\mathcal{G}$ yields a $2$fi-optimal row-column design if and only if
no column of $\mathcal{G}_c$ or $\mathcal{G}_r$ is a zero vector.
\end{proposition}

\section{Construction of 2fi-optimal  $s^n$ full factorial row-column designs}

In this section, we will present a theoretical construction of  2fi-optimal  $s^n$ full factorial row-column designs
for any odd prime level $s$,
as stated in Theorem 4.
Notice that if $p+q=n$, an $s^n$ full factorial row-column design arranged in an $s^p\times s^q$ array is obtained
from a $(p+q)\times (p+q)$ array generator matrix.
Since an array generator matrix has full row rank,
it has full column rank as well, so all columns are linearly independent.

To describe our constructions, we fix the following notation which will be used in the rest of this paper:

$I_v$: the $v\times v$ identity matrix;

$J_v$: the $v\times v$ matrix with all retries $1$;

$\boldsymbol{0}_{v\times w}$: the $v\times w$ matrix with all entries $0$;

$\boldsymbol{1}_v$: the $v\times 1$ matrix with all entices $1$;

$\boldsymbol{e_i}$: an column vector with $i$-th entry 1 and other entries 0;

$H_v$: the $v\times v$ matrix $(h_{ij})_{v\times v}$ with $h_{ij}=1$ if $i+j=v+1$ and $h_{ij}=0$ otherwise.

\begin{theorem} \label{full}
Let $s$ be an odd prime. Suppose $n=p+q$. Then
 $\mathcal{G}$ yields a $2$fi-optimal $s^{n}$ full factorial row-column design if
\begin{itemize}
\item[(a)] for $1=p\leq q$,
$$\mathcal{G}=\left(
\begin{array}{c}
\mathcal{G}_c\\  \hline
\mathcal{G}_r\\
\end{array}
\right)=\left(
\begin{array}{c:c}
1 & \boldsymbol{1}_q^T\\ \hline
\boldsymbol{1}_q & I_q+J_q \\
\end{array}
\right);$$

\item[(b)] for $2=p\leq q$,
$$\mathcal{G}=\left(
\begin{array}{c}
\mathcal{G}_c\\  \hline
\mathcal{G}_r\\
\end{array}
\right)=\left(
\begin{array}{c:c:c}
I_2 & M & X_{*}\\ \hline
I_2&M+I_2 &\boldsymbol{0}_{2\times (q-2)}\\\hdashline
\boldsymbol{0}_{(q-2)\times 2}&\boldsymbol{0}_{(q-2)\times 2}&I_{q-2}\\
\end{array}
\right),$$
where $M=\left(\begin{array}{cc}
1 & 1\\
2 &1 \\
\end{array}
\right)$,
$X_{*}$ is a $2\times (q-2)$ matrix which is selected so that each point of $\rm{PG(1,s)}$ appears in
the $2\times n$ matrix
$\begin{array}{c;{3pt/1pt}c;{3pt/1pt}c}
(I_2 &   M &   X_{*})
\end{array}$,
$\alpha$ or $\alpha+1$ times, where $\alpha=\lfloor n/(s+1)\rfloor$, and the columns of $X_{*}$ are not zero vectors$;$
\item[(c)] for $3\leq p\leq q$, $$\mathcal{G}=\left(
\begin{array}{c}
\mathcal{G}_c\\  \hline
\mathcal{G}_r\\
\end{array}
\right)=\left(
\begin{array}{c:c:c}
I_p & J_p+H_p & Y_{*}\\ \hline
H_p&J_p+2I_p &\boldsymbol{0}_{p\times(q-p)}\\\hdashline
\boldsymbol{0}_{(q-p)\times p}&\boldsymbol{0}_{(q-p)\times p}&I_{q-p}\\
\end{array}
\right),$$
where $Y_{*}$ is a $p\times (q-p)$ matrix which is selected so that each point of $\rm{PG(p-1,s)}$ appears in
the $p\times n$ matrix
$\begin{array}{c;{3pt/1pt}c;{3pt/1pt}c}(I_p & J_p+H_p &   Y_{*})
\end{array}$,
$\alpha$ or $\alpha+1$ times, where $\alpha=\lfloor n(s-1)/(s^p-1)\rfloor$, and the columns of $Y_{*}$ are not zero vectors.
\end{itemize}
\end{theorem}

\proof
(a) Matrix $\mathcal{G}$ has row reduced echelon form  $\left(
\begin{array}{c:c}
1 & \boldsymbol{1}_q^T\\ \hline
\boldsymbol{0}_{q\times 1} & I_q \\
\end{array}
\right)$ working modulo $s$, so it is full rank.
By Proposition \ref{p=1}, case (a) holds.

(b) Reducing $\mathcal{G}$ to row echelon form, working modulo $s$, gives $$\left(
\begin{array}{c:c:c}
I_2 & M & X_{*}\\ \hline
\boldsymbol{0}_{2\times 2}&I_2 &-X_{*}\\ \hdashline
\boldsymbol{0}_{(q-2)\times 2}&\boldsymbol{0}_{(q-2)\times 2}&I_{q-2}\\
\end{array}
\right),$$  which confirms that $\mathcal{G}$ is full rank.
Condition (1) of Proposition \ref{conditions} is satisfied because no column in $\mathcal{G}_c$ or $\mathcal{G}_r$
is a zero vector.
Since the columns of $\begin{array}{c;{3pt/1pt}c}
(I_2 &M) \end{array}$ belong to $4$ different points of $\rm{PG(1,s)}$,  the existence of $X_{*}$ is guaranteed and thus
Condition (2) is satisfied.
All columns of $\mathcal{G}_r$ belong to different points of $\rm{PG(q-1,s)}$, so Condition (3) also holds.

(c) Reducing $\mathcal{G}$ to row echelon form, working modulo $s$, gives $$\left(
\begin{array}{c:c:c}
I_p & J_p+H_p & Y_{*}\\ \hline
\boldsymbol{0}_{p\times p}&I_p &-H_pY_{*}\\ \hdashline
\boldsymbol{0}_{(q-p)\times p}&\boldsymbol{0}_{(q-p)\times p}&I_{q-p}\\
\end{array}
\right),$$  which confirms that $\mathcal{G}$ is full rank.
Condition (1) of Proposition \ref{conditions} is satisfied because no column in $\mathcal{G}_c$ or $\mathcal{G}_r$
is a zero vector.
The columns of $\begin{array}{c;{3pt/1pt}c}
(I_p &J_p+H_p) \end{array}$ are  different points of $\rm{PG(p-1,s)}$ when $p\geq 3$, so the existence of $Y_{*}$ is proved and thus
Condition (2) is satisfied.
Condition (3) also holds because all columns of $\mathcal{G}_r$ belong to different points of $\rm{PG(q-1,s)}$.\qed

Some specific constructions of Theorem \ref {full} are given in the Appendix.

\section{Construction of 2fi-optimal  $s^{n-1}$ fractional factorial row-column designs}

Full factorial experiments require a prohibitively large number of runs when the number of factors is large.
One alternative is to use the fractional factorial designs which have fewer runs and
allow estimation of all main effects and some two-factor interactions.
For $p+q<n$, the row-column design constructed by an array generator matrix of order $(p+q)\times n$
is a $s^{(p+q)-n}$ fraction of the full $s^n$ factorial.
In this section, we focus on constructions of $s^{n-1}$ fractional factorial row-column designs for any prime $s$,
which is divided into two subsections, one is for two-level, the other is for odd prime level, both of which
cover all combinations of parameters $(s,p,q,n)$ under the restriction $n=p+q+1$.

\subsection{2-level constructions}

For $p=1$ and $q\leq2$, there is at least one confounded main effect,
so no construction exists in the case. Thus if $p=1$, we assume $q\geq 3$.

\subsubsection{ Constructions for $p=1$}

\begin{theorem}\label{2level p=1}
Suppose $n=p+q+1$.
 For $p=1, q\geq 3$, $\mathcal{G}$ yields a $2$fi-optimal $2^{n-1}$ fractional factorial row-column design if
$$\mathcal{G}=\left(
\begin{array}{c}
\mathcal{G}_c\\  \hline
\mathcal{G}_r\\
\end{array}
\right)=\left(
\begin{array}{c:c:c}
1 & \boldsymbol{1}_q^T &1\\ \hline
\boldsymbol{1}_q& I_q+J_q &\boldsymbol{e}_q\\
\end{array}
\right).$$
\end{theorem}
\proof
By reduction of matrix $\mathcal{G}$ to row echelon form, working modulo $2$, we get matrix
 $$\mathcal{G}_0=\left(
\begin{array}{c:c:c}
1 & \boldsymbol{1}_q^T &1\\ \hline
\boldsymbol{0}_{q\times 1}& I_q&\boldsymbol{e}_q-\boldsymbol{1}_q \\
\end{array}
\right).$$
 We claim that any three columns of $\mathcal{G}_0$ are linearly independent.
It's easy to see that the first $q+1$ columns of $\mathcal{G}_0$ are linearly independent.
Suppose to the contrary that there exist  $3$ column of $\mathcal{G}_0$ being  linearly dependent.
Then the only possible case is that the last column  of $\mathcal{G}_0$ is a linear combination
of $2$ columns from the first $q+1$ columns of $\mathcal{G}_0$. Thus $\boldsymbol{e}_q-\boldsymbol{1}_q$
is a linear combination of $2$ columns of
$\begin{array}{c;{3pt/1pt}c}(\boldsymbol{0}_{q\times 1} & I_q)\end{array}$. As $q\geq 3$,
this  implies that $q=3$ and $\mathcal{G}_0$ turns into
$$\left(
\begin{array}{c:ccc:c}
1 & 1&1&1 &1\\ \hline
0&1&0&0&1\\
0&0&1&0&1\\
0&0&0&1&0\\
\end{array}
\right),$$
which forces $(1,1,1,0)^T=(1,1,0,0)^T+(1,0,1,0)^T$ working modulo $2$, a contradiction.
Hence the claim is proved. Therefore, any three columns of $\mathcal{G}$ are also linearly independent.
Since no column of $\mathcal{G}_c$ or $\mathcal{G}_r$ is a zero vector, by Proposition \ref{p=1},
the design generated by $\mathcal{G}$ is $2$fi-optimal.\qed

\subsubsection{Constructions for $p=2$}
\begin{theorem} \label{2level p=2}
Suppose $n=p+q+1$. $\mathcal{G}$ yields a $2$fi-optimal $2^{n-1}$ fractional factorial row-column design if
\begin{itemize}
\item[(a)] for $p=2$, $q=2$,
$$\mathcal{G}=\left(
\begin{array}{c}
\mathcal{G}_c\\  \hline
\mathcal{G}_r\\
\end{array}
\right)=\left(
\begin{array}{cc:cc:c}
1&1&0&1&1\\
0&1&1&1&1\\ \hline
1&1&1&0&1\\
1&0&1&1&1\\
\end{array}
\right);$$

\item[(b)] for $p=2$, $q=3$,
$$\mathcal{G}=\left(
\begin{array}{c}
\mathcal{G}_c\\  \hline
\mathcal{G}_r\\
\end{array}
\right)=\left(
\begin{array}{c:c:c}
I_2 & F & \boldsymbol{1}_2\\ \hline
E& I_3+EF &\boldsymbol{e}_3\\
\end{array}
\right),$$
where $E=\left(\begin{array}{cc}
1 & 0\\
0 &1 \\
1 &1 \\
\end{array}
\right)$,
$F=\left(\begin{array}{ccc}
1 & 0 &1\\
1 &1 &0\\
\end{array}
\right)$;

\item[(c)] for $p=2$, $q= 4$,
$$\mathcal{G}=\left(
\begin{array}{c}
\mathcal{G}_c\\  \hline
\mathcal{G}_r\\
\end{array}
\right)=\left(
\begin{array}{c:c:c:c}
I_2 & F & \boldsymbol{1}_2 &\boldsymbol{e}_1\\ \hline
E& I_3+EF &\boldsymbol{e}_1 &\boldsymbol{0}_{3\times 1}\\
\hdashline
\boldsymbol{0}_{1\times 2} &\boldsymbol{0}_{1\times 3} & 1 &1\\
\end{array}
\right),$$
where $E$ and $F$ are as in case $(b)$;

\item[(d)] for $p=2$, $q\geq 5$,
$$\mathcal{G}=\left(
\begin{array}{c}
\mathcal{G}_c\\  \hline
\mathcal{G}_r\\
\end{array}
\right)=\left(
\begin{array}{c:c:c:c}
I_2 & F & \boldsymbol{1}_2  & X_{*}\\ \hline
E & I_3+EF &\boldsymbol{1}_3 & \boldsymbol{0}_{3\times(q-3)} \\ \hdashline
\boldsymbol{0}_{(q-3)\times 2} &\boldsymbol{0}_{(q-3)\times 3} & \boldsymbol{1}_{q-3} &I_{q-3}\\
\end{array}
\right),$$
where $E, F$ are as in case $(b)$, $X_{*}$ is a $2\times (q-3)$ matrix which is selected so that each point of $\rm{PG(1,2)}$ appears in
the $2\times n$ matrix
$\begin{array}{c;{3pt/1pt}c;{3pt/1pt}c;{3pt/1pt}c}(I_2 & F& \boldsymbol{1}_2 &  X_{*})
\end{array}$,
$\alpha$ or $\alpha+1$ times, where $\alpha=\lfloor n/3\rfloor$, and the columns of $X_{*}$ are not zero vectors.

\end{itemize}
\end{theorem}

\proof
(a) By adding the first row to the third row of $\mathcal{G}$, and by adding the sum of the first two rows to the fourth row of $\mathcal{G}$,
working modulo $2$, we have the row reduced echelon form
$$\mathcal{G}_0=\left(\begin{array}{cc:cc:c}
1&1&0&1&1\\
0&1&1&1&1\\ \hline
0&0&1&1&0\\
0&0&0&1&1\\
\end{array}\right),$$
which is full row rank and so $\mathcal{G}$ is also full row rank. This implies that the first $4$ columns of $\mathcal{G}$ are linearly independent.
Observe that the last column of $\mathcal{G}$ is the sum of the first $4$ columns, so it cannot be linearly
represented by any three of them.
Therefore  any four columns of $\mathcal{G}$ are linearly independent.
Denote the $5$ columns of $\mathcal{G}$ by factors $A, B,C,D,E$, respectively. Then by Lemma \ref{identify} the two-factor interactions $AB$,
$AD$, $BC$, $CD$ are unconfounded. By Zhou and Zhou (2023, Theorem 1), we have the upper bound $\phi (2,2,2,5)=8$, and
the $2$fi-efficiency for $p=q=s=2$, $n=5$ is $0.5$. Matrix $\mathcal{G}$
enables  all main effects  and $4$ two-factor interactions to be unconfounded and  thus is  $2$fi-optimal.

(b) Reducing $\mathcal{G}$ to row echelon form, working modulo $2$,  gives
$$\mathcal{G}_0=\left(
\begin{array}{c:c:c}
I_2 & F & \boldsymbol{1}_2\\ \hline
\boldsymbol{0}_{3\times 2}& I_3 &\boldsymbol{1}_3\\
\end{array}
\right).$$
We claim that any four columns of $\mathcal{G}_0$ are linearly independent.
It's easy to see that the first $5$ columns of $\mathcal{G}_0$ are linearly independent.
Suppose to the contrary that there exist $4$ vectors of $\mathcal{G}_0$ being  linearly dependent.
The only possibility is that the last column of $\mathcal{G}_0$ is a linear combination
of $3$ columns from the first $5$ columns of $\mathcal{G}_0$, which implies that
$(1,1,1,1,1)^T=(1,1,1,0,0)^T+(0,1,0,1,0)^T+(1,0,0,0,1)^T$, a contradiction.
Thus the claim is true and hence any four columns of $\mathcal{G}$ are also linearly independent.
Denote the $6$ columns of $\mathcal{G}$ by factors $A, B,C,D,E, F$, respectively. Then by Lemma \ref{identify},
the two-factor interactions $AB$, $AC$, $AD$, $AF$, $BC$, $BE$, $BF$, $CD$, $CE$, $DE$, $EF$
are unconfounded. By Zhou and Zhou (2023, Theorem 1), the upper bound $\phi(s,p,q,n)=\phi (2,2,3,6)=12$
and the $2$fi-efficiency for $p=2, q=3, s=2, n=6$ is $0.9167$. Matrix $\mathcal{G}$
enables  all main effects  and $11$ two-factor interactions to be unconfounded and  thus is  $2$fi-optimal.

(c) Reducing $\mathcal{G}$ to row echelon form, working modulo $2$, gives
$$\mathcal{G}_0=\left(
\begin{array}{c:c:c:c}
I_2 & F & \boldsymbol{1}_2 &\boldsymbol{e}_1\\ \hline
\boldsymbol{0}_{3\times 2}& I_3 &\boldsymbol{e}_2 &\boldsymbol{e}_{1}+\boldsymbol{e}_3\\
\hdashline
\boldsymbol{0}_{1\times 2} &\boldsymbol{0}_{1\times 3} & 1 &1\\
\end{array}
\right).$$
We claim that any four columns of $\mathcal{G}_0$ are linearly independent.
It's easy to see that the first $6$ columns of $\mathcal{G}_0$ are linearly independent.
Suppose to the contrary that there exist $4$ vectors of $\mathcal{G}_0$ being  linearly dependent.
The only possibility is that the last column of $\mathcal{G}_0$ is a linear combination
of $3$ columns from the first $6$ columns of $\mathcal{G}_0$. By taking the last $4$ entries of each vector,  we see this happens only if
$(1,0,1,1)^T=(1,0,0,0)^T+(0,0,1,0)^T+(0,1,0,1)^T$, a contradiction.
Thus the claim is true and hence any four columns of $\mathcal{G}$ are also linearly independent.
Clearly, Condition (1) of Proposition \ref{conditions} is satisfied.
The number of columns of $\mathcal{G}_c$ belonging to each point of $\rm{PG(1,2)}$ is $2$ or $3$,  thus satisfying Condition (2).
Since the columns of

$$\begin{array}{c;{3pt/1pt}c}(E & I_3+EF )\end{array}
=\left(\begin{array}{cc:ccc}
1&0&0&0&1\\
0&1&1&0&0\\
1&1&0&1&0\\
\end{array}\right)$$
are all different,  the columns of $\mathcal{G}_r$ are all different, satisfying Condition (3).

(d) Reducing $\mathcal{G}$ to row echelon form, working modulo $2$, yields
$$\mathcal{G}_0
=\left(
\begin{array}{c:c:c:c}
I_2 & F & \boldsymbol{1}_2  & X_{*}\\ \hline
\boldsymbol{0}_{3\times 2} & I_3 &\boldsymbol{e}_3 & EX_{*} \\ \hdashline
\boldsymbol{0}_{(q-3)\times 2} &\boldsymbol{0}_{(q-3)\times 3} & \boldsymbol{1}_{q-3} &I_{q-3}\\
\end{array}
\right),$$which implies that the first $6$ columns of $\mathcal{G}$ are linearly independent.
We claim that any four columns of $\mathcal{G}$ are linearly independent.
Let $\mathcal{L}$ denote the set of the first $5$ columns of $\mathcal{G}$,
and let $\mathcal{R}$ denote the set of the last $q-2$ columns of $\mathcal{G}$.
It's easy to see that all $5$ vectors in $\mathcal{L}$ are linearly independent.
Suppose to the contrary that there exist $4$ vectors of $\mathcal{L}\cup\mathcal{R}$ being  linearly dependent,
among which there are $t$ vectors from $\mathcal{L}$, and $4-t$ vectors from $\mathcal{R}$, where $0\leq t\leq 3$.
By checking the bottom right-hand $q\times {(q-2)}$ submatrix  of $\mathcal{G}$,  we  derive that the last $q-2$ columns of $\mathcal{G}$
are linearly independent, so $t\neq 0$.
Checking the last $q-3$ rows of $\mathcal{G}$ gives $t\neq 2, 3$.  If $t=1$ then $q=5$, which implies
that the sum of the last $3$ columns of $\mathcal{G}$ is a vector in $\mathcal{L}$. This further means
$ \begin{array}{c;{3pt/1pt}c}
 (E & I_3+EF)
 \end{array}$ has a column equal to $(1,1,1)^T$, which is a contradiction
by checking the $3$-th to $5$-th entries of the columns of $\mathcal{G}$.
Clearly,  Condition (1) of Proposition \ref{conditions} is satisfied. Since the number of columns of
$\begin{array}{c;{3pt/1pt}c;{3pt/1pt}c}(I_2 & F& \boldsymbol{1}_2)
\end{array}$
belonging to each point of $\rm{PG(1,2)}$ is $2$,  the existence of $X_{*}$ is guaranteed,
and so Condition (2) is satisfied.
Since the columns of
$\begin{array}{c;{3pt/1pt}c}(E & I_3+EF )\end{array}$
are all different,  the columns of $\mathcal{G}_r$ are all different, satisfying Condition (3). \qed

\subsubsection{Constructions for $ p\geq 3$}

\begin{lemma}\label{K_p}
For $p\geq 3$, let $K_p=(k_{i,j})$ be the $p\times p$ matrix with $k_{i,j}=1$ if $i\leq j$ and $k_{i,j}=0$ if $i>j$,
and let $L_p=(l_{i,j})=I_p+K_p(I_p+J_p)$. Then the columns of the $p\times (2p+1)$ matrix $\begin{array}{c;{3pt/1pt}c;{3pt/1pt}c}
(K_p& L_p &K_p\boldsymbol{1}_p+\boldsymbol{1}_p)\end{array}$ are  all
 distinct and nonzero vectors.
\end{lemma}

\proof
By Godolphin (2019, Lemma 1), the columns of  $\begin{array}{c;{3pt/1pt}c}
(K_p& L_p )\end{array}$ are all distinct and nonzero vectors.
Thus we only need to prove that the vector $K_p\boldsymbol{1}_p+\boldsymbol{1}_p$ is not a column of $K_p$ or $L_p$.
Note that the last entry of $K_p\boldsymbol{1}_p+\boldsymbol{1}_p$ is $0$, while the last
row of $L_p$ is $\boldsymbol{1}_p^T$, so $K_p\boldsymbol{1}_p+\boldsymbol{1}_p$ is not a column of $L_p$.
Taking the last $3$ entries of $K_p\boldsymbol{1}_p+\boldsymbol{1}_p$ gives $(0,1,0)^T$, while taking the last
$3$ entries of each column of $K_p$ gives $(0,0,0)^T$, $(1,0,0)^T$, $(1,1,0)^T$ or $(1,1,1)^T$,
so $K_p\boldsymbol{1}_p+\boldsymbol{1}_p$ is not a column of $K_p$.\qed


\begin{theorem}\label{2level p=3}
Suppose $n=p+q+1$.
For $3\leq p\leq q$, $\mathcal{G}$ yields a $2$fi-optimal $2^{n-1}$ fractional factorial row-column design if

$$\mathcal{G}=\left(
\begin{array}{c}
\mathcal{G}_c\\  \hline
\mathcal{G}_r\\
\end{array}
\right)=\left(
\begin{array}{c:c:c:c}
I_p & I_p+J_p  & \boldsymbol{1}_p & X_{*}\\ \hline
K_p & L_p &K_p\boldsymbol{1}_p+\boldsymbol{1}_p & \boldsymbol{0}_{p\times(q-p)} \\ \hdashline
\boldsymbol{0}_{(q-p)\times p} &\boldsymbol{0}_{(q-p)\times p} & \boldsymbol{0}_{(q-p)\times 1} &I_{q-p}\\
\end{array}
\right),$$
where
$X_{*}$ is a $p\times (q-p)$ matrix which is selected so that each point of $\rm{PG(p-1,2)}$ appears in
the $p\times n$ matrix
$\begin{array}{c;{3pt/1pt}c;{3pt/1pt}c;{3pt/1pt}c}(I_p & I_p+J_p& \boldsymbol{1}_p &  X_{*})
\end{array}$,
$\alpha$ or $\alpha+1$ times, where $\alpha=\lfloor n/(2^p-1)\rfloor$, and the columns of $X_{*}$ are not zero vectors.
\end{theorem}

\proof
Reducing $\mathcal{G}$ to row echelon form, working modulo $2$,  gives
$$\mathcal{G}_0=\left(
\begin{array}{c:c:c:c}
I_p & I_p+J_p  & \boldsymbol{1}_p & X_{*}\\ \hline
\boldsymbol{0}_{p\times p} & I_p & \boldsymbol{1}_p & K_pX_{*}\\ \hdashline
\boldsymbol{0}_{(q-p)\times p} &\boldsymbol{0}_{(q-p)\times p} & \boldsymbol{0}_{(q-p)\times 1} &I_{q-p}\\
\end{array}
\right),$$which implies that the first $2p$ columns of $\mathcal{G}_0$ are linearly independent.
We claim that any four columns of $\mathcal{G}_0$ are linearly independent.
First assume $p=q$. If there exist 4 linearly dependent columns of $\mathcal{G}_0$,
then the $2p+1$-th column is a linear combination of three columns from the first $2p$ columns of $\mathcal{G}_0$.
This implies that $p=3$ and so $\boldsymbol{1}_3$ is the sum of 3 columns of $I_3+J_3$,
a contradiction. Thus assume $q>p$ in the following.
Let $\mathcal{L}$ denote the set of the first $2p$ columns of $\mathcal{G}_0$,
and let $\mathcal{R}$ denote the set of the last $q-p+1$ columns of $\mathcal{G}_0$.
Suppose to the contrary that there exist $4$ vectors of $\mathcal{L}\cup\mathcal{R}$ being  linearly dependent,
among which there are $t$ vectors from $\mathcal{L}$, and $4-t$ vectors from $\mathcal{R}$, where $0\leq t\leq 3$.
By checking the right bottom $q\times {(q-p+1)}$ submatrix  of $\mathcal{G}_0$,  we  derive that the last $q-p+1$ columns of $\mathcal{G}_0$
are linearly independent, so $t\neq 0$.
Checking the last $q-p$ rows of $\mathcal{G}_0$, one easily excludes the cases $t=1$ and $2$.
If $t=3$, then the only possibility is that the $(2p+1)$-th column of $\mathcal{G}_0$ is a linear combination of $3$ vectors of $\mathcal{L}$.
Suppose $p=3$.   It follows that the $7$-th column is the sum of
the $4$-th, $5$-th and $6$-th columns of $\mathcal{G}_0$,
so the sum of columns of $I_3+J_3$ is equal to $\boldsymbol{1}_3$, a contradiction.
Checking the $p$-th to $(2p)$-th entries of the columns of $\mathcal{G}_0$ excludes the case $p>3$.

Therefore, the claim is true and hence any four columns of $\mathcal{G}$ are also linearly independent.
Clearly, no column of $\mathcal{G}_c$ or $\mathcal{G}_r$ is a zero vector.
The columns of $\begin{array}{c;{3pt/1pt}c;{3pt/1pt}c}(I_p & I_p+J_p& \boldsymbol{1}_p )
\end{array}$ are all distinct, so $X_{*}$ exists.
From Lemma \ref{K_p}, we derive that the columns of $\mathcal{G}_r$ are all different.
Therefore matrix $\mathcal{G}$ satisfies conditions of Proposition \ref{conditions}, so the design generated by $\mathcal{G}$
is $2$fi-optimal and the $2$fi-efficiency is $1$. \qed

Some construction examples of this subsection are postponed in the Appendix.

\subsection{Odd prime level constructions}
Suppose $s$ is an odd prime.
If $p=q=1$, then  there exists at least one confounded main effect,
so no construction exists in the case.

\subsubsection{Constructions for $p=1$}

\begin{theorem}\label{fractional p=1}
 Let $s$ be an odd prime and suppose $n=p+q+1$. For $p=1, q\geq 2$, $\mathcal{G}$ yields a $2$fi-optimal $s^{n-1}$ fractional factorial row-column design if
$$\mathcal{G}=\left(
\begin{array}{c}
\mathcal{G}_c\\  \hline
\mathcal{G}_r\\
\end{array}
\right)=\left(
\begin{array}{c:c:c}
1 & \boldsymbol{1}_q^T &2\\ \hline
\boldsymbol{1}_q& I_q+J_q &\boldsymbol{1}_q\\
\end{array}
\right).$$
\end{theorem}

\proof
By reduction of matrix $\mathcal{G}$ to row echelon form, working modulo $s$, we get matrix
 $$\mathcal{G}_0=\left(
\begin{array}{c:c:c}
1 & \boldsymbol{1}_q^T &2\\ \hline
\boldsymbol{0}_{q\times 1}& I_q&-\boldsymbol{1}_q \\
\end{array}
\right).$$
 We claim that any three columns of $\mathcal{G}_0$ are linearly independent.
It's easy to see that the first $q+1$ columns of $\mathcal{G}_0$ are linearly independent.
Suppose to the contrary that there exist $3$ columns of $\mathcal{G}_0$ being  linearly dependent.
Then the only possible case is that the last column  of $\mathcal{G}_0$ is a linear combination
of $2$ columns from the first $q+1$ columns of $\mathcal{G}_0$. Thus $-\boldsymbol{1}_q$
is a linear combination of $2$ columns of
$\begin{array}{c;{3pt/1pt}c}(\boldsymbol{0}_{q\times 1} & I_q)\end{array}$. As $q\geq 2$,
this  implies that $q=2$ and $\mathcal{G}_0$ turns into
$$\left(
\begin{array}{c:cc:c}
1 & 1&1 &2\\ \hline
0&1&0&-1\\
0&0&1&-1\\
\end{array}
\right),$$
which forces $(2,-1,-1)^T=-(1,1,0)^T-(1,0,1)^T$  and thus $2=-2 \pmod{s}$ , a contradiction as $s$ is an odd prime.
Hence the claim is proved. Therefore, any three columns of $\mathcal{G}$ are also linearly independent.
Since no columns of $\mathcal{G}_c$ or $\mathcal{G}_r$ is a zero vector, by Proposition \ref{p=1},
the design generated by $\mathcal{G}$ is $2$fi-optimal. \qed

\subsubsection{Constructions for $p=2$}
\begin{theorem}\label{fractional p=2}
 Let $s$ be an odd prime and suppose $n=p+q+1$. For $q\geq p=2$, $\mathcal{G}$ yields a $2$fi-optimal $s^{n-1}$ fractional factorial row-column design if
\begin{itemize}
\item[(a)] for $p=2, q=2$,
$$\mathcal{G}=\left(
\begin{array}{c}
\mathcal{G}_c\\  \hline
\mathcal{G}_r\\
\end{array}
\right)=\left(
\begin{array}{c:c:c}
I_2 & M & \boldsymbol{a}\\ \hline
I_2 & M+I_2 & \boldsymbol{b}\\
\end{array}
\right),$$
where $\boldsymbol{a}=(1,2)^T$, $\boldsymbol{b}=(2,0)^T$, $M=\left(
\begin{array}{cc}
1 &1\\
2&1\\
\end{array}
\right)$ if $s=3$, and
 $\boldsymbol{a}=(1,1)^T$, $\boldsymbol{b}=(2,2)^T$, $M=\left(
\begin{array}{cc}
1 &1\\
3&2\\
\end{array}
\right)$ otherwise;

\item[(b)] for $p=2, q=3$,
$$\mathcal{G}=\left(
\begin{array}{c}
\mathcal{G}_c\\  \hline
\mathcal{G}_r\\
\end{array}
\right)=\left(
\begin{array}{c:c:c:c}
I_2 & M & \boldsymbol{a}& \boldsymbol{b}\\ \hline
I_2&M+I_2 &\boldsymbol{e}_1&\boldsymbol{0}_{2\times 1}\\\hdashline
\boldsymbol{0}_{1\times 2}&\boldsymbol{0}_{1\times 2}&1&1\\
\end{array}
\right),$$
where $M=\left(\begin{array}{cc}
1 & 1\\
2 &1 \\
\end{array}
\right)$,
$\boldsymbol{a}=(1,s-2)^T$ and $\boldsymbol{b}=(1,s-1)^T$;

\item[(c)] for $p=2, q= 4$,
$$\mathcal{G}=\left(
\begin{array}{c}
\mathcal{G}_c\\  \hline
\mathcal{G}_r\\
\end{array}
\right)=\left(
\begin{array}{c:c:c:c:c}
I_2 & M & \boldsymbol{a} &\boldsymbol{b}&\boldsymbol{c}\\ \hline
I_2& M+I_2 &\boldsymbol{e}_2 &\boldsymbol{0}_{2\times 1}&\boldsymbol{0}_{2\times 1}\\
\hdashline
\boldsymbol{0}_{2\times2} &\boldsymbol{0}_{2\times 2} & \boldsymbol{1}_2 &\boldsymbol{e}_1+2\boldsymbol{e}_2&\boldsymbol{e}_2\\
\end{array}
\right),$$
where $M$ is as in case (b), $\boldsymbol{a}=(1,s-1)^T$, $\boldsymbol{b}=(1,s-2)^T$ and $\boldsymbol{c}=(1,s-3)^T$;

\item[(d)] for $p=2, q\geq 5$,
$$\mathcal{G}=\left(
\begin{array}{c}
\mathcal{G}_c\\  \hline
\mathcal{G}_r\\
\end{array}
\right)=\left(
\begin{array}{c:c:c:c}
I_2 & M & \boldsymbol{a}  & X_{*}\\ \hline
I_2 & M+I_2 &\boldsymbol{a} & \boldsymbol{0}_{2\times(q-2)} \\ \hdashline
\boldsymbol{0}_{(q-2)\times 2} &\boldsymbol{0}_{(q-2)\times 2} & \boldsymbol{1}_{q-2} &I_{q-2}\\
\end{array}
\right),$$
where $M$ is as in case (b), $\boldsymbol{a}=(1,s-1)^T$, $X_{*}$ is a $2\times (q-2)$ matrix which is selected so that each point of $\rm{PG(1,s)}$ appears in
the $2\times n$ matrix
$\begin{array}{c;{3pt/1pt}c;{3pt/1pt}c;{3pt/1pt}c}(I_2 & M& \boldsymbol{a} &   X_{*})\end{array}$,
$\alpha$ or $\alpha+1$ times, where $\alpha=\lfloor n/(s+1)\rfloor$, and the columns of $X_{*}$ are not zero vectors.

\end{itemize}
\end{theorem}

\proof
(a) If $s=3$, reducing $\mathcal{G}$ to row echelon form, working modulo $s$, gives
$$\mathcal{G}_0=\left(
\begin{array}{cc:cc:c}
1& 0& 1&1 & 1\\
0&1&2&1&2\\ \hline
0&0&1&0&1\\
0&0&0&1&1\\
\end{array}
\right).$$
It is clear that the first $4$ columns of $\mathcal{G}_0$ are linearly independent.
If the last column of $\mathcal{G}_0$ is a linear combination of other $3$ columns of $\mathcal{G}_0$,
then the only possibility is that
$(1,2,1,1)^T=(1,2,1,0)^T+(1,1,0,1)^T+\boldsymbol{u}$, where $\boldsymbol{u}$ is
a multiple of the first or second column of $\mathcal{G}_0$, which is not possible as $\boldsymbol{u}=(-1,-1,0,0)^T$.
Therefore, any $4$ columns of $\mathcal{G}_0$ and hence of $\mathcal{G}$ are linearly independent.
Denote the $5$ columns of $\mathcal{G}$ by factors $A, B,C,D,E$, respectively. Then by Lemma \ref{identify}, all
  $10$ two-factor interactions are unconfounded  except $A\times E$ and $C\times E$.
By Zhou and Zhou (2022, Theorem 1), we have the upper bound $\phi(s,p,q,n)=\phi (3,2,2,5)=9$, and
 by the $17$-th row of Table 1 (Zhou and Zhou, 2022), the $2$fi-efficiency for $p=2, q=2, s=3, n=5$ is $0.8889$. Matrix $\mathcal{G}$
enables  all main effects  and $8$ two-factor interactions to be unconfounded and  thus is  $2$fi-optimal.

If $s>3$, reducing $\mathcal{G}$ to row echelon form, working modulo $s$, gives
$$\mathcal{G}_0=\left(
\begin{array}{cc:cc:c}
1& 0& 1&1 & 1\\
0&1&3&2&1\\ \hline
0&0&1&0&1\\
0&0&0&1&1\\
\end{array}
\right).$$
It is clear that the first $4$ columns of $\mathcal{G}_0$ are linearly independent.
If the last column of $\mathcal{G}_0$ is a linear combination of other $3$ columns of $\mathcal{G}_0$,
then the only possibility is that
$(1,1,1,1)^T=(1,3,1,0)^T+(1,2,0,1)^T+\boldsymbol{u}$, where $\boldsymbol{u}$ is
a multiple of the first or second column of $\mathcal{G}_0$, which is not possible as $\boldsymbol{u}=(-1,-4,0,0)^T$.
Therefore, any $4$ columns of $\mathcal{G}_0$ and hence of $\mathcal{G}$ are linearly independent.
None of the columns of $\mathcal{G}_c$ or $\mathcal{G}_r$ is a zero vector,  so Condition (1) of Proposition \ref{conditions}
is satisfied. All $5$ columns of $\mathcal{G}_c$
belong to  distinct points of $\rm{PG(1,s)}$, and all $5$ columns of $\mathcal{G}_r$
also belong to distinct points of $\rm{PG(1,s)}$, so Conditions (2) and (3) are satisfied.
Therefore, the design generated by $\mathcal{G}$ is $2$fi-optimal and the $2$fi-efficiency is $1$.

(b) Reducing $\mathcal{G}$ to row echelon form, working modulo $s$,  gives
$$\mathcal{G}_0=\left(
\begin{array}{c:c:c:c}
I_2 & M & \boldsymbol{a}& \boldsymbol{b}\\ \hline
\boldsymbol{0}_{2\times 2}&I_2 &\boldsymbol{e}_1-\boldsymbol{a}&-\boldsymbol{b}\\ \hdashline
\boldsymbol{0}_{1\times 2}&\boldsymbol{0}_{1\times 2}&1&1\\
\end{array}
\right).$$
It is clear that the first $5$ columns of $\mathcal{G}_0$ are linearly independent.
If the last column of $\mathcal{G}_0$ is a linear combination of other $3$ columns of $\mathcal{G}_0$,
then the only possibility is that
$(1,-1,-1,1,1)^T=(1,-2,0,2,1)^T+\lambda\boldsymbol{u}+\mu\boldsymbol{v}$, where $\lambda,\mu\in \rm{GF(s)}$,  $\boldsymbol{u}$ is
 the third column of $\mathcal{G}_0$,
and $\boldsymbol{v}$ is the fourth column of $\mathcal{G}_0$. Checking the third and fourth rows of $\mathcal{G}_0$
gives $\lambda=\mu=-1$, which further implies $(1,-1,-1,1,1)^T=(1,-2,0,2,1)^T-(1,2,1,0,0)^T-(1,1,0,1,0)^T$,  a contradiction.
Therefore any $4$ columns of $\mathcal{G}_0$ and hence of $\mathcal{G}$
are linearly independent.
None of the columns of $\mathcal{G}_c$ or $\mathcal{G}_r$ is a zero vector,  so Condition (1) of Proposition \ref{conditions}
is satisfied.  Each point of $\rm{PG(1,s)}$ appears in the columns of $\mathcal{G}_c$ $\alpha$ or $\alpha+1$ times,
where $\alpha=0$ if $s>5$, and $\alpha=1$ if $s=3$ or $5$.
All $6$ columns of $\mathcal{G}_r$
belong to distinct points of $\rm{PG(2,s)}$, so Conditions (2) and (3) are satisfied.
Therefore, the design generated by $\mathcal{G}$ is $2$fi-optimal and the $2$fi-efficiency is $1$.

(c) Reducing $\mathcal{G}$ to row echelon form, working modulo $s$,  gives
$$\mathcal{G}_0=\left(
\begin{array}{c:c:c:c:c}
I_2 & M & \boldsymbol{a} &\boldsymbol{b}&\boldsymbol{c}\\ \hline
\boldsymbol{0}_{2\times 2}& I_2 &\boldsymbol{e}_2 -\boldsymbol{a}&-\boldsymbol{b}&-\boldsymbol{c}\\
\hdashline
\boldsymbol{0}_{2\times2} &\boldsymbol{0}_{2\times 2} & \boldsymbol{1}_2 &\boldsymbol{e}_1+2\boldsymbol{e}_2&\boldsymbol{e}_2\\
\end{array}
\right).$$
It is clear that the first $6$ columns of $\mathcal{G}_0$ are linearly independent.
Suppose that the last column of $\mathcal{G}_0$ is a linear combination of other $3$ columns of $\mathcal{G}_0$.
Observe that $\boldsymbol{e}_2$ has a unique linear representation by $\boldsymbol{1}_2$ and $\boldsymbol{e}_1+2\boldsymbol{e}_2$,
which is given by $\boldsymbol{e}_2=-\boldsymbol{1}_2+(\boldsymbol{e}_1+2\boldsymbol{e}_2)$.
This forces that the last column of $\mathcal{G}_0$ equals $-\boldsymbol{u}+\boldsymbol{v}+\boldsymbol{w}$,
where $\boldsymbol{u}$ is the $5$-th column, $\boldsymbol{v}$ is the $6$-th column and $\boldsymbol{w}$
is a multiple of one of  the first $4$ columns of $\mathcal{G}_0$. This implies $\boldsymbol{w}=(1,-2,-1,3,0,0)^T$,
which is impossible.
 Therefore any $4$ columns of $\mathcal{G}_0$ and hence of $\mathcal{G}$
are linearly independent.
None of the columns of $\mathcal{G}_c$ or $\mathcal{G}_r$ is a zero vector,  so Condition (1) of Proposition \ref{conditions}
is satisfied.  Each point of $\rm{PG(1,s)}$ appears in the columns of $\mathcal{G}_c$ $\alpha$ or $\alpha+1$ times,
where $\alpha=0$ if $s>5$, and $\alpha=1$ if $s=3$ or $5$.
All $7$ columns of $\mathcal{G}_r$
belong to  distinct points of $\rm{PG(3,s)}$, so Conditions (2) and (3) are satisfied.
Therefore, the design generated by $\mathcal{G}$ is $2$fi-optimal and the $2$fi-efficiency is $1$.

(d) Reducing $\mathcal{G}$ to row echelon form, working modulo $s$,  gives
$$\mathcal{G}_0=\left(
\begin{array}{c:c:c:c}
I_2 & M & \boldsymbol{a}  & X_{*}\\ \hline
\boldsymbol{0}_{2\times 2} &I_2 &\boldsymbol{0}_{2\times 1} & -X_{*} \\ \hdashline
\boldsymbol{0}_{(q-2)\times 2} &\boldsymbol{0}_{(q-2)\times 2} & \boldsymbol{1}_{q-2} &I_{q-2}\\
\end{array}
\right),$$
which implies that the first $4$ columns of $\mathcal{G}$ are linearly independent.
We claim that any four columns of $\mathcal{G}$ are linearly independent.
Let $\mathcal{L}$ denote the set of the first $4$ columns of $\mathcal{G}$,
and let $\mathcal{R}$ denote the set of the last $q-1$ columns of $\mathcal{G}$.
Suppose to the contrary that there exist $4$ vectors of $\mathcal{L}\cup\mathcal{R}$ being  linearly dependent,
among which there are $t$ vectors from $\mathcal{L}$, and $4-t$ vectors from $\mathcal{R}$, where $0\leq t\leq 3$.
By checking the bottom right-hand $q\times {(q-1)}$ submatrix  of $\mathcal{G}$,  we  derive that the last $q-1$ columns of $\mathcal{G}$
are linearly independent, so $t\neq 0$.
 Checking the last $q-2$ rows of $\mathcal{G}$ gives $t\neq 1, 2, 3$. Thus the claim is proved.
Clearly,  Condition (1) of Proposition \ref{conditions} is satisfied.
Each point of $\rm{PG(1,s)}$ appears in the columns of
$\begin{array}{c;{3pt/1pt}c;{3pt/1pt}c}(I_2 & M & \boldsymbol{a})
\end{array}$ $\gamma$ or $\gamma +1$ times, where $\gamma=1$ if $s=3$, and $\gamma=0$ if $s>3$.
So the existence of $X_{*}$ is guaranteed,  and thus Condition (2) is satisfied.
Since the columns of
$\begin{array}{c;{3pt/1pt}c}(I_2 & M+I_2 )\end{array}$
are not multiples of each other,  the columns of $\mathcal{G}_r$ belong to distinct points of $\rm{PG(q-1,s)}$, satisfying Condition (3).
Therefore  the design generated by $\mathcal{G}$
is $2$fi-optimal and the $2$fi-efficiency is $1$. \qed

\subsubsection{Constructions for $p\geq 3$}

\begin{theorem}\label{fractional p=3}
Let $s$ be an odd prime and suppose $n=p+q+1$.  For $3=p\leq q$, $\mathcal{G}$ yields a $2$fi-optimal $s^{n-1}$ fractional factorial row-column design if

$$\mathcal{G}=\left(
\begin{array}{c}
\mathcal{G}_c\\  \hline
\mathcal{G}_r\\
\end{array}
\right)=\left(
\begin{array}{c:c:c:c}
I_p & J_p+H_p  & \boldsymbol{a}& X_{*}\\ \hline
H_p & J_p+2I_p &\boldsymbol{b} & \boldsymbol{0}_{p\times(q-p)} \\ \hdashline
\boldsymbol{0}_{(q-p)\times p} &\boldsymbol{0}_{(q-p)\times p} & \boldsymbol{0}_{(q-p)\times 1} &I_{q-p}\\
\end{array}
\right),$$
where $\boldsymbol{a}=\boldsymbol{1}_p$, $\boldsymbol{b}=2\boldsymbol{a}+\boldsymbol{e}_1-\boldsymbol{e}_p$,
$X_{*}$ is a $p\times (q-p)$ matrix which is selected so that each point of $\rm{PG(p-1,s)}$ appears in
the $p\times n$ matrix $\begin{array}{c;{3pt/1pt}c;{3pt/1pt}c;{3pt/1pt}c}(I_p & J_p+H_p& \boldsymbol{a} &  X_{*})\end{array}$, $\alpha$ or $\alpha+1$ times, where $\alpha=\lfloor n(s-1)/(s^p-1)\rfloor$, and the columns of $X_{*}$ are not zero vectors.

\end{theorem}

\proof
Reducing $\mathcal{G}$ to row echelon form, working modulo $s$,  gives
$$\mathcal{G}_0=\left(
\begin{array}{c:c:c:c}
I_p & J_p+H_p  & \boldsymbol{a}& X_{*}\\ \hline
\boldsymbol{0}_{p\times p} & I_p &\boldsymbol{b}-H_p\boldsymbol{a} & -H_pX_{*} \\ \hdashline
\boldsymbol{0}_{(q-p)\times p} &\boldsymbol{0}_{(q-p)\times p} & \boldsymbol{0}_{(q-p)\times 1} &I_{q-p}\\
\end{array}
\right),$$
which implies that all  columns of $\mathcal{G}_0$ excluding the $(2p+1)$-th column are linearly independent.
Suppose that there exist $4$ columns of $\mathcal{G}_0$ being  linearly dependent. The only possibility is that
the $(2p+1)$-th column is a linear combination of other $3$ columns from the first $2p$ columns of $\mathcal{G}_0$.
This is impossible for $p\geq 5$ as $\boldsymbol{b}-H_p\boldsymbol{a}$ has $p-1$ nonzero entries.
If $p=3$, then $\boldsymbol{b}-H_p\boldsymbol{a}=(2,1,0)^T$. This forces
$(1,1,1)^T=2(1,1,2)^T+(1,2,1)^T+\boldsymbol{u}$, where $\boldsymbol{u}$ is a multiple of
$e_1, e_2$ or $e_3$.
However, this equation gives $\boldsymbol{u}=(s-2,s-3,s-4)$, a contradiction.
If $p=4$, then $\boldsymbol{b}-H_p\boldsymbol{a}=(2,1,1,0)^T$. This forces
$(1,1,1,1)^T=2(1,1,1,2)^T+(1,1,2,1)^T+(1,2,1,1)^T$, which is  impossible working modulo $s$.
Thus any $4$ columns of $\mathcal{G}_0$ and hence of $\mathcal{G}$ are linearly independent.
Clearly,  Condition (1) of Proposition \ref{conditions} is satisfied.
The columns of
$\begin{array}{c;{3pt/1pt}c;{3pt/1pt}c}(I_p & J_p+H_p & \boldsymbol{a})
\end{array}$ are not multiples of each other when $p\geq 3$,  and this confirms the existence of $X_{*}$.
  Thus Condition (2) is satisfied.
Since the columns of
$\begin{array}{c;{3pt/1pt}c;{3pt/1pt}c}(H_p & J_p+2I_p &\boldsymbol{b} )\end{array}$
are not multiples of each other,  the columns of $\mathcal{G}_r$ belong to distinct points of $\rm{PG(q-1,s)}$, satisfying Condition (3).
Therefore the design generated by $\mathcal{G}$
is $2$fi-optimal and the $2$fi-efficiency is $1$. \qed

Theorems 9 --11 provide many 2fi-optimal constructions, see some examples in the Appendix.

\section{Discussion}

This work focuses on estimation of main effects and two-factor interactions in row-column designs in a higher prime level. The array generator matrix
provides immediate and illuminating information on confounded treatment effects by Lemma \ref{identify}, and
the use of concepts from finite projective geometry facilitates the statements of sufficient conditions of 2fi-optimal designs.
The 2fi-optimality is an important criterion for evaluating row-column designs.
This paper gives systematic constructions of 2fi-optimal $s^n$ full factorial and $s^{n-1}$ fractional factorial row-column designs for any parameter
combination with $s$ prime. Further research could focus on systematic constructions of 2fi-optimal $s^{n-k}$ fractional
factorial row-column designs
for $k\geq 2$.

\begin{appendices}
\section*{Appendix}

\subsection*{Examples of Theorem \ref{full}}
\begin{itemize}
\item [(a)] $s\geq3, p=1,q=2, n=3$,
$$\mathcal{G}=\left(
\begin{array}{c:cc}
1 & 1&1\\ \hline
1&2&1\\
1&1&2
\end{array}
\right).$$

\item [(b)] $s=3, p=2,q=3,n=5$,
$$\mathcal{G}=\left(
\begin{array}{cc:cc:c}
1 &0& 1&1&1\\
0&1&2&1&1\\ \hline
1&0&2&1&0\\
0&1&2&2&0\\  \hdashline
0&0&0&0&1
\end{array}
\right).$$



\item [(c)] $s=5, p=3,q=4, n=7$,
$$\mathcal{G}=\left(
\begin{array}{ccc:ccc:c}
1 &0&0& 1&1&2&1\\
0&1&0&1&2&1&1\\
0&0&1&2&1&1&1\\ \hline
0&0&1&3&1&1&0\\
0&1&0&1&3&1&0\\
1&0&0&1&1&3&0\\ \hdashline
0&0&0&0&0&0&1
\end{array}
\right).$$
\end{itemize}

\subsection*{An example of Theorem \ref{2level p=1}} \

$s=2, p=1, q=3, n=5$,
$$\mathcal{G}=\left(
\begin{array}{c:ccc:c}
1 &1&1& 1&1\\ \hline
1&0&1&1&0\\
1&1&0&1&0\\
1&1&1&0&1
\end{array}
\right).$$

\subsection*{Examples of Theorem \ref{2level p=2}}
\begin{itemize}
\item[(a)] $s=2, p=2,q=3, n=6$,
$$\mathcal{G}=\left(
\begin{array}{cc:ccc:c}
1 &0&1& 0&1&1\\
0&1&1&1&0&1\\ \hline
1&0&0&0&1&0\\
0&1&1&0&0&0\\
1&1&0&1&0&1
\end{array}
\right).$$
\item [(b)] $s=2, p=2,q=4, n=7$,
$$\mathcal{G}=\left(
\begin{array}{cc:ccc:c:c}
1 &0&1& 0&1&1&1\\
0&1&1&1&0&1&0\\ \hline
1&0&0&0&1&1&0\\
0&1&1&0&0&0&0\\
1&1&0&1&0&0&0\\ \hdashline
0&0&0&0&0&1&1
\end{array}
\right).$$

\item [(c)] $s=2, p=2,q=5, n=8$,
$$\mathcal{G}=\left(
\begin{array}{cc:ccc:c:cc}
1 &0&1& 0&1&1&1&0\\
0&1&1&1&0&1&0&1\\ \hline
1&0&0&0&1&1&0&0\\
0&1&1&0&0&1&0&0\\
1&1&0&1&0&1&0&0\\ \hdashline
0&0&0&0&0&1&1&0\\
0&0&0&0&0&1&0&1
\end{array}
\right).$$

\end{itemize}

\subsection*{An example of Theorem \ref{2level p=3}}\

$s=2,p=3,q=4, n=8$,
$$\mathcal{G}=\left(
\begin{array}{ccc:ccc:c:c}
1 &0&0& 0&1&1&1&1\\
0&1&0&1&0&1&1&0\\
0&0&1&1&1&0&1&0\\ \hline
1&1&1&1&0&0&0&0\\
0&1&1&0&0&1&1&0\\
0&0&1&1&1&1&0&0\\ \hdashline
0&0&0&0&0&0&0&1
\end{array}
\right).$$

\subsection*{An example of Theorem \ref{fractional p=1}}\

$s\geq3, p=1, q=2, n=4$,
$$\mathcal{G}=\left(
\begin{array}{c:cc:c}
1 &1&1&2\\ \hline
1&2&1&1\\
1&1&2&1
\end{array}
\right).$$

\subsection*{Examples of Theorem \ref{fractional p=2}}
\begin{itemize}

\item[(a)] $s=5, p=2, q=2, n=5$,
$$\mathcal{G}=\left(
\begin{array}{cc:cc:c}
1 &0&1& 1&1\\
0&1&3&2&1\\ \hline
1&0&2&1&2\\
0&1&3&3&2
\end{array}
\right).$$

\item [(b)] $s=3, p=2,q=3, n=6$,
$$\mathcal{G}=\left(
\begin{array}{cc:cc:c:c}
1 &0&1& 1&1&1\\
0&1&2&1&1&2\\ \hline
1&0&2&1&1&0\\
0&1&2&2&0&0\\ \hdashline
0&0&0&0&1&1
\end{array}
\right).$$



\item [(c)] $s=5, p=2,q=4, n=7$,
$$\mathcal{G}=\left(
\begin{array}{cc:cc:c:c:c}
1 &0&1& 1&1&1&1\\
0&1&2&1&4&3&2\\ \hline
1&0&2&1&0&0&0\\
0&1&2&2&1&0&0\\\hdashline
0&0&0&0&1&1&0\\
0&0&0&0&1&2&1
\end{array}
\right).$$


\item [(d)] $s=5, p=2,q=5, n=8$,
$$\mathcal{G}=\left(
\begin{array}{cc:cc:c:ccc}
1 &0&1& 1&1&1&0&1\\
0&1&2&1&4&0&1&3\\ \hline
1&0&2&1&1&0&0&0\\
0&1&2&2&4&0&0&0\\  \hdashline
0&0&0&0&1&1&0&0\\
0&0&0&0&1&0&1&0\\
0&0&0&0&1&0&0&1
\end{array}
\right).$$
\end{itemize}

\subsection*{Examples of Theorem \ref{fractional p=3}}\
\begin{itemize}
\item[(a)]
$s=3,p=3,q=4, n=8,$
$$\mathcal{G}=\left(
\begin{array}{ccc:ccc:c:c}
1 &0&0& 1&1&2&1&1\\
0&1&0&1&2&1&1&1\\
0&0&1&2&1&1&1&0\\ \hline
0&0&1&0&1&1&0&0\\
0&1&0&1&0&1&2&0\\
1&0&0&1&1&0&1&0\\ \hdashline
0&0&0&0&0&0&0&1
\end{array}
\right).$$

\item[(b)]
$s=5,p=3,q=4, n=8,$
$$\mathcal{G}=\left(
\begin{array}{ccc:ccc:c:c}
1 &0&0& 1&1&2&1&1\\
0&1&0&1&2&1&1&1\\
0&0&1&2&1&1&1&0\\ \hline
0&0&1&3&1&1&3&0\\
0&1&0&1&3&1&2&0\\
1&0&0&1&1&3&1&0\\ \hdashline
0&0&0&0&0&0&0&1
\end{array}
\right).$$
\end{itemize}
\end{appendices}

\section*{Funding}
This paper was partially supported by the National Natural Science Foundation of China (11961076, 12271471).

\end{document}